\documentclass[a4paper]{article}
\usepackage{cite}
\usepackage[all]{xy}
\usepackage[latin1]{inputenc}      
\usepackage[dvips]{graphics,graphicx}
\usepackage{appendix}
\usepackage{amsfonts,amssymb,amsmath,color, mathrsfs, amstext}
\usepackage{amsbsy, amsopn, amscd, amsxtra, amsthm,authblk}
\usepackage{enumerate,algorithmicx,algorithm}
\usepackage{algpseudocode}
\usepackage{upref, xcolor}
\usepackage{geometry,wrapfig}
\usepackage{amsthm,amsmath,amssymb}
\usepackage{relsize}
\usepackage{mathrsfs}
\usepackage{bm}
\usepackage{ulem}
\usepackage{exscale}
\usepackage{tabularx}
\usepackage{threeparttable,multirow,dcolumn,booktabs}
\usepackage{algorithmicx,algorithm}
\usepackage{tabularx}
\usepackage{caption}
\usepackage[displaymath]{lineno}
\usepackage {subfigure} 
\usepackage{float}
\usepackage{bm}
\usepackage{caption}
\usepackage{yhmath}
\usepackage{booktabs}
\usepackage{multirow}
\usepackage{makecell}
\usepackage[colorlinks,
            linkcolor=red,
            anchorcolor=red,
            citecolor=red
            ]{hyperref}
\usepackage{cleveref}
\numberwithin{equation}{section}

\numberwithin{equation}{section}

\newcommand{\bs}[1]{\boldsymbol{#1}}
\newcommand{\bmu}{\bs{\mu}}

\newcommand{\calN}{\mathcal{N}}

\numberwithin{equation}{section}
\usepackage{epstopdf}
\usepackage{pdfpages}

\graphicspath{{./Figs/}}

\begin{document}

\title{Derivative-informed Graph Convolutional Autoencoder with Phase Classification for the Lifshitz-Petrich Model}

\author{
Yanlai Chen\footnote{Department of Mathematics, University of Massachusetts Dartmouth, North Dartmouth, MA 02747. Email: {\tt{yanlai.chen@umassd.edu}}. Y. Chen is partially supported by National Science Foundation grant DMS-2208277.}, \quad  
Yajie Ji\footnote{Department of Statistics and Data Science, Yale University, New Haven, CT 06511. Email: {\tt yajie.ji@yale.edu}. School of Mathematical Sciences, Shanghai Jiao Tong University, Shanghai 200240, China. Email: {\tt jiyajie595@sjtu.edu.cn}. Y. Ji acknowledges the support from the NSFC (No. 124B2023).},\quad
Zhenli Xu\footnote{School of Mathematical Sciences, CMA-Shanghai and MOE-LSC, Shanghai Jiao Tong University, Shanghai 200240, China. Email: {\tt xuzl@sjtu.edu.cn}. Z. Xu acknowledges the support from the NSFC (grant Nos. 12325113 and 12426304) and the SJTU Kunpeng \& Ascend Center of Excellence.}
}

\date{}
\makeatother

\maketitle

\begin{abstract}
The Lifshitz-Petrich (LP) model is a classical model for describing complex spatial patterns such as quasicrystals and multiphase structures. Solving and classifying the solutions of the LP model is challenging due to the presence of high-order gradient terms and the long-range orientational order characteristic of the quasicrystals. To address these challenges, we propose a Derivative-informed Graph Convolutional Autoencoder (DiGCA) {to classify the} multi-component multi-state  solutions of the LP model. The classifier consists of two stages. In the offline stage, the DiGCA phase classifier innovatively incorporates both solutions and their derivatives for training a graph convolutional autoencoder which effectively captures intricate spatial dependencies while significantly reducing the dimensionality of the solution space. In the online phase, the framework employs a neural network classifier to efficiently categorize encoded solutions into distinct phase diagrams. The numerical results demonstrate that the DiGCA phase classifier accurately solves the LP model, classifies its solutions, and rapidly generates detailed phase diagrams in a robust manner, offering significant improvements in both efficiency and accuracy over traditional methods.

\end{abstract}

\section{Introduction}

The Lifshitz-Petrich (LP) model \cite{lifshitz1997theoretical} provides a crucial theoretical framework for exploring complex spatial patterns, such as quasicrystals and other aperiodic structures. It introduces a free energy functional incorporating terms that account for both the interaction energy and the bulk energy with two characteristic scales, such that the corresponding partial differential equation (PDE)  models multiple symmetries within the system. The interaction energy term incorporates high-order differential operators to induce the formation of ordered structures\cite{yin2021transition} such as crystals, quasicrystals, and multiphase crystals, while the remaining components correspond to bulk energy contributions, typically expressed through polynomial or logarithmic functions. The LP model is instrumental in understanding self-assembly processes in various physical and material systems, including condensed matter physics, metallurgy, and soft matter\cite{hayashida2007polymeric,percec2009self,talapin2009quasicrystalline,zeng2004supramolecular,zhang2012dodecagonal}. Its primary significance lies in its ability to describe the formation and stability of intricate patterns through a set of high-order parametric partial differential equations (pPDEs). There are two characteristic parameters, with one representing the temperature and the other delineating the asymmetry of the order parameter.

Traditional numerical methods for solving the LP model, including finite difference, finite element, and spectral methods, face considerable challenges due to the high dimensionality and complexity of the solution space. Among these approaches, the crystalline approximation method (CAM)\cite{lifshitz1997theoretical,zhang2008efficient} and the projection method (PM)\cite{jiang2014numerical} are particularly popular due to their effectiveness in handling the intricate patterns and dynamics described by the LP Model. The CAM approximates quasicrystals using a periodic structure which is called Simultaneous Diophantine Approximation in number theory, but as the accuracy of the approximation increases, the size of the computational domain grows rapidly, leading to a prohibitive computational cost and making the method impractical for high-accuracy simulations\cite{jiang2014numerical}. On the other hand, the PM leverages the fact that the reciprocal vectors of a quasicrystal in a lower-dimensional space can be approximated by linear combinations of basic reciprocal vectors in a higher-dimensional space, effectively rendering the periodic property in a higher-dimensional space. However, the computational burden increases exponentially with dimensionality, eventually rendering the problem numerically intractable\cite{jiang2014numerical,ji2024mcms}. These methods often necessitate substantial computational resources and advanced algorithms to ensure both accuracy and stability \cite{gao2025nc}. In general, while full order models based on traditional discretizations are known for their high accuracy, they require extensive computational resources, which often become prohibitive for parametrized problems. 

To address this challenge, surrogate models, such as reduced basis models \cite{Quarteroni2015} are developed as more efficient emulators to computationally expensive solvers. These surrogate models strike a balance between accuracy and efficiency, allowing faster simulations while maintaining acceptable levels of accuracy\cite{franco2023deep}. Among these models, machine learning methods, particularly deep neural networks and deep neural operator (such as DeepONet\cite{lu2021deepOnet} and Fourier Neural Operator (FNO)\cite{li2021fourier}), have demonstrated significant progress in solving PDEs and pPDEs. 
Based on whether access to the full-order model is required, these reduction methods can be broadly categorized into intrusive and non-intrusive approaches. Intrusive methods, such as Proper Orthogonal Decomposition\cite{chatterjee2000podintroduction} and the Reduced Basis Method (RBM)\cite{Quarteroni2015}, construct reduced basis spaces by selecting representative solutions or parameters using techniques such as SVD or greedy algorithms. These methods then project the high-dimensional discretized equations onto a low-dimensional function space through approaches like Galerkin projection which is often used to ensure that the error is minimized within this subspace. Intrusive methods are well-suited for systems governed by explicit and relatively simple physical equations but face limitations when applied to complex nonlinear systems or high-dimensional parameterized problems. In contrast, non-intrusive methods directly construct reduced-order models using solutions from the full-order models without requiring access to the original model. Data-driven non-intrusive approaches include neural operator methods such as DeepONet and FNO, which learn mappings from parameters to solutions, and neural network-based nonlinear dimensionality reduction techniques such as classic auto-encoders\cite{hinton2006reducing}, convolutional auto-encoders\cite{masci2011stacked}, variational auto-encoders\cite{pu2016variational}. Furthermore, hybrid approaches such as POD-NN \cite{hesthaven2018nonNNRBM} and POD-DL-ROM \cite{fresca2021comprehensive,fresca2022pod} combine traditional POD methods for spatial dimensionality reduction with neural networks to map the reduced parameter space to the solution space.

In this work, we follow Pichi {\it et al.}\cite{pichi2024graph} and introduce the Derivative-informed Graph Convolutional Autoencoder (DiGCA) phase classifier for the LP model, a novel two-step method designed to efficiently generate phase diagrams. The graph convolutional autoencoders \cite{pichi2024graph} 
exploit the flexibility of graph-based representations to capture intricate spatial dependencies and reduce the dimensionality of high-dimensional data. This capability renders graph convolutional autoencoders particularly suitable for addressing the challenges posed by the LP model. The proposed DiGCA consists of two subnetworks, each comprising an offline and an online stage. In the offline stage of the first subnetwork, we utilize both the solutions of the LP model and their derivatives to train a multi-component multi-state graph convolutional autoencoders for each stable state, referred to as MCMS-DiGCA. The MCMS-DiGCA captures detailed spatial dependencies and trains tailored networks for different stable states, identifying the corresponding mappings for each. In the online stage, MCMS-DiGCA can efficiently provide solutions and high-order derivative information corresponding to given parameters. In the second subnetwork, the DiGCA method employs a phase classifier to classify the encoded solutions from MCMS-DiGCA into distinct phase diagrams. Through supervised training with labeled data in the offline stage, the online stage of this second subnetwork leverages the outputs of MCMS-DiGCA to successfully predict phase diagrams with an accuracy exceeding $98\%$. Compared to traditional model order reduction methods like the
multi-component multi-state reduced basis method (MCMS-RBM), the DiGCA with phase classification method is highly computationally efficient, achieving a two-order-of-magnitude improvement in computational speed, which marks a substantial advancement in applying machine learning techniques to high-order PDEs, offering an efficient framework for exploring complex physical systems with multiple states. Moreover, we numerically demonstrate that the method is robust by showing that the phase diagram accuracy suffers little to no degradation with up to 10\% of additive white noise.

We note that the accurate computation of high-order derivatives and the efficient training of models to capture these derivatives are non-trivial tasks. Traditional neural network methods rely on automatic differentiation and backpropagation to handle high-order derivatives, typically requiring complex computations and substantial memory resources, significantly increasing computational complexity and storage demands. Chen {\it et al.}\cite{lyu2022mim} proposed a deep mixed residual method (MIM), which introduces a novel approach to solving high-order PDEs by reformulating them into first-order systems, inspired by classical methods such as local discontinuous Galerkin and mixed finite element methods. MIM employs the least squares residual of the first-order system as the loss function, reducing the computational burden of high-order derivatives while offering flexibility through various loss functions and neural network architectures. Our novel two-step DiGCA phase classifier augments the input to the vanilla Convolutional Authoencoder by a judiciously selected set of derivatives.

The rest of this paper is organized as follows. In Section \ref{sec:bgd}, we introduce the parametrized LP model, the projection method, and the pseudospectral method used to obtain the high-fidelity solution. The key components of the DiGCA and phase classification, including the online and offline procedures, are introduced in Section \ref{sec:GAPC}. We then present numerical results in Section \ref{sec:numerics} to demonstrate the efficiency, accuracy, and robustness of the proposed DiGCA with phase classification. Finally, concluding remarks are drawn in Section \ref{sec:conclusion}.

\section{Lifshitz-Petrich model}
\label{sec:bgd}
In this section, we introduce the  Lifshitz-Petrich model and the process preparing the dataset based on the projection method.

The Lifshitz-Petrich (LP) model extends two fundamental pattern-formation models: the Swift-Hohenberg model\cite{swift1977hydrodynamic}  with free energy functional
\[
\mathcal{F}_{\text{SH}}(\phi;c,\varepsilon,\alpha)=\int_V d r\left\{\frac{c}{2}\left|\left(\nabla^2+1^2\right) \phi\right|^2-\frac{\varepsilon}{2} \phi^2+\frac{1}{4} \phi^4\right\}
\]
widely used in materials science, and the Landau-Brazovskii (LB) model \cite{brazovskii1975phase}
\[\mathcal{F}_{\text{LB}}(\phi;c,\varepsilon,\alpha)=\int_V d r\left\{\frac{c}{2}\left|\left(\nabla^2+1^2\right) \phi\right|^2-\frac{\varepsilon}{2} \phi^2-\frac{\alpha}{3} \phi^3+\frac{1}{4} \phi^4\right\}
\]
which describes polymeric systems. It extends the Landau-Brazovskii model by incorporating two characteristic wavelength scales, allowing for more complex pattern formations. The model incorporates several parameters, one representing temperature and the other characterizing the asymmetry of the order parameter. The scalar order parameter $\phi(\bm{r})$ describes how perfectly the molecules are aligned. For each parameter configuration, the LP model yields a unique steady-state solution as well as multiple metastable solutions. Notably, the steady state corresponds to the global minimum of the free energy functional that is divided into two contributions,
\begin{equation}
\mathcal{F}(\phi; c, q, \varepsilon,\alpha)= E_1(\phi; c, q) + E_2(\phi; \varepsilon, \alpha),
\label{eq:phi_energy}
\end{equation}
where $\bm r \in \mathcal{R}^d$ with $d=2$, $V$ is the system volume, and
\begin{equation}
\label{eq:GH-defined}
\begin{aligned}
& E_1(\phi; c, q) =\int_V \frac{c}{2}\lvert\mathcal{G}(\phi; q)\rvert^2d \bm{r}, 
& E_2(\phi; \varepsilon, \alpha) = \int_V \mathcal{H}(\phi; \varepsilon, \alpha)d \bm{r}\\
& \mathcal{G}(\phi; q) = \left(\nabla^2+1^2\right)\left(\nabla^2+q^2\right) \phi, 
& \mathcal{H}(\phi; \varepsilon, \alpha) = -\frac{\varepsilon}{2} \phi^2-\frac{\alpha}{3} \phi^3+\frac{1}{4} \phi^4.
\end{aligned}
\end{equation}

We see that the free energy functional in Eq.~\eqref{eq:phi_energy} includes two main components: the nonlocal term $\mathcal{H}$ expressed using polynomial or logarithmic functions, and the local term $\mathcal{G}$ which employs high-order differential operators to generate structures such as crystals, quasicrystals, and disorder state.

A critical feature of the LP model is the use of an energy penalty parameter $c=1$ to ensure that the principal reciprocal vectors of structures are located on $\lvert \bm{k}\rvert=1$ and $\lvert \bm{k}\rvert=q = 2\cos(\pi/12)$, with $q$ being an irrational number depending on the symmetry. The reduced temperature $\varepsilon$ and the phenomenological parameter $\alpha$ are the key variables of interest in this study for phase diagram construction.
 
For a given parameter $\bmu:=[\varepsilon,\alpha]$, the candidate stable states correspond to the local minima of the free energy functional. These states are determined by solving the Euler-Lagrange equation associated with the functional,
\begin{equation}
\frac{\delta \mathcal{F}}{\delta \phi}= 0.
\end{equation}
This is an eighth-order nonlinear partial differential equation.
To solve it, the gradient flow method can be employed \cite{gradientflowshen2019,ju2022generalized}, such that the order parameter obeys
\begin{equation}
\frac{\partial \phi}{\partial t}=-\frac{\delta \mathcal{F}}{\delta \phi}=-\left[c\left(\nabla^2+1\right)^2\left(\nabla^2+q^2\right)^2 \phi-\varepsilon \phi\right] -\frac{\alpha}{3} \phi^3+\frac{1}{4} \phi^4.
\label{eq:gradient_flow_form}
\end{equation}
By the following implicit-explicit discretization scheme for the time, one has,
\begin{equation}
\frac{\phi^{n+1}-\phi^n}{\Delta t}=-\left[c\left(\nabla^2+1\right)^2\left(\nabla^2+q^2\right)^2 \phi^{n+1}-\varepsilon \phi^n\right] -\frac{\alpha}{3} (\phi^n)^3+\frac{1}{4} (\phi^n)^4.
\label{eq:semi_implicit}
\end{equation}

\begin{figure}[htbp]
\centering
\includegraphics[width=1\textwidth]{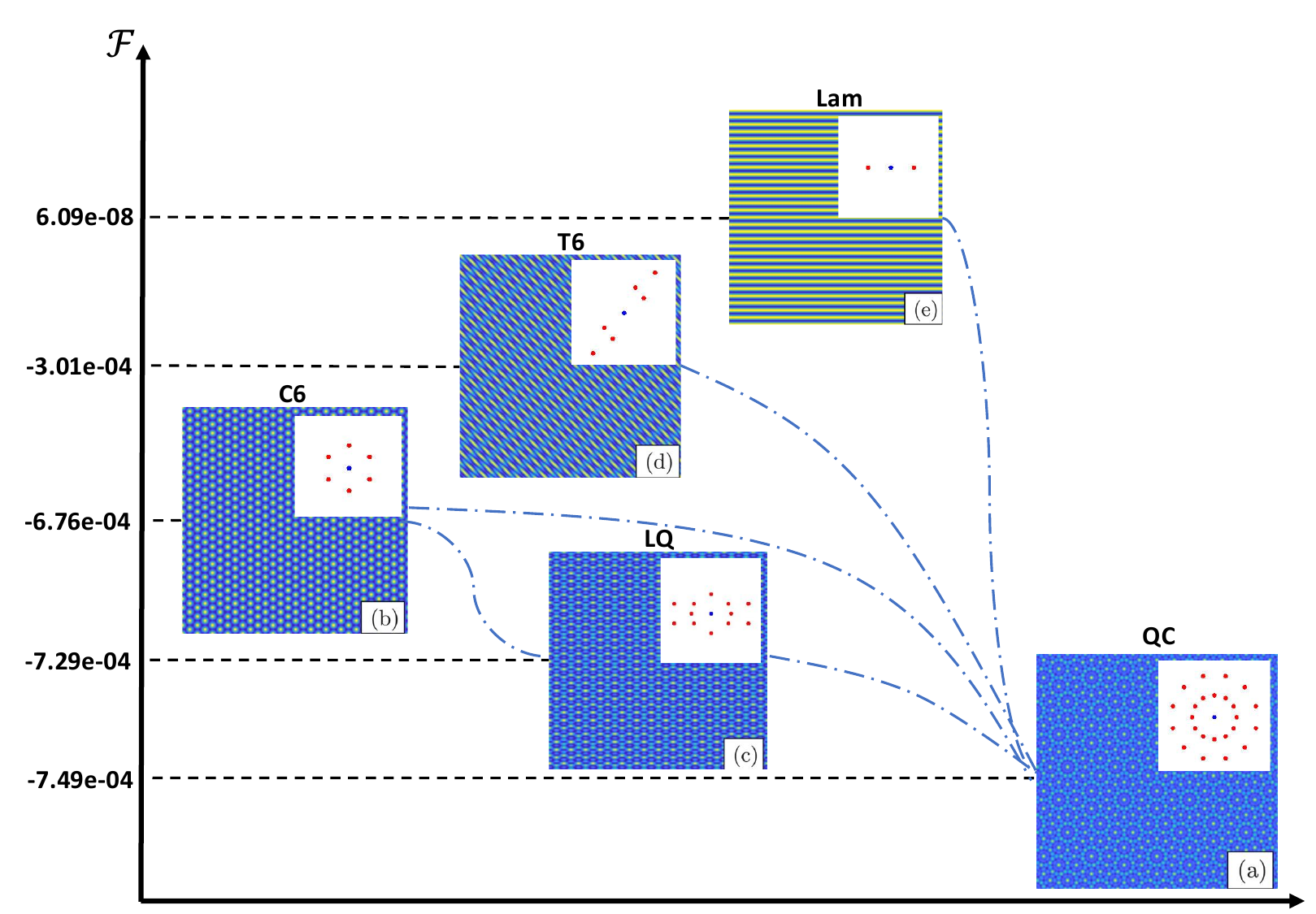}
\caption{Order parameters with the corresponding prominent diffraction patterns in the reciprocal space for QC (a), C6(b), LQ(c), T6(d), and Lam(e) states. The top figure also shows the stable (QC) and metastable solutions(C6, LQ, T6, Lam) corresponding to parameter $(\varepsilon,\alpha) = (5\times 10^{-6},\sqrt2/2)$, their respective energies, and possible energy transfer pathways (denoted by blue lines).}
\label{fig:five_states}
\end{figure} 

By varying these parameters $\bmu$, the order parameter $\phi$ exhibits a diverse range of equilibrium ordered phases, resulting in a complex phase behavior and leading to quasicrystals (QC), periodic structures or liquid phase (Lq). If one takes $q=2\cos (\pi/12)$ in the free energy functional $\mathcal{F}(\phi(\bm{r}); c, q, \varepsilon,\alpha)$, the LP model can explain the 12-fold symmetry QC excited by dual-frequency filtering in the Faraday experiment that provides a macroscopic description of the system\cite{lifshitz1997theoretical}.  This 12-fold QC in two dimensions is the stable state with the minimum energy. In addition to the stable state QC, the LP model has been shown to have different metastable states, including the 6-fold crystalline state (C6), the lamellar quasicrystalline state (LQ), the transformed 6-fold crystalline state (T6), and the Lamella state (Lam) \cite{yin2021transition}. Figure \ref{fig:five_states} illustrates the structures of these stable states with parameter $(c,\varepsilon,\alpha) = (1,5\times 10^{-6},\sqrt2/2)$ in both real and reciprocal spaces. Among these five states, QC was first discovered by Shechtman $\textit{et al.}$ in the 1980s with a 5-fold rotational symmetries in a rapidly cooling Al-Mn alloy\cite{shechtman1984metallic}. After that, more different structures with 5-, 6-, 8-, 12-, and 20-fold symmetries have emerged in various metallic alloys.

In general, numerical computation of quasiperiodic problems employs either the CAM or PM, both of which transform quasiperiodic issues into periodic ones, for which efficient and advanced algorithms are well-established. Specifically, the pseudospectral method achieves efficiency by evaluating the gradient terms in the Fourier space and the nonlinear terms in the physical space. The PM \cite{jiang2014numerical} initially constructs periodic structures in a higher-dimensional reciprocal space, which are subsequently projected into lower-dimensional space through a projection matrix $\mathcal{S}$. According to the LP model's dimensional constraint, we set $n=4$ and $d=2$ as the working parameters. Notably, in the 12-fold symmetric system, particular reciprocal vectors resist representation as integer-coefficient linear combinations of the two-dimensional basis vectors $\bm{p}_1^*=(0,1)$ and $\bm{p}_4^*=(1,0)$, i.e., 
\begin{equation}
    \bm{k} \neq k_1\bm{p}_1^* + k_4\bm{p}_4^* \quad \text{for any}\ k_1,k_4 \in \mathbb{Z}.
\end{equation}
This prevents direct application of Fourier analysis. To resolve this limitation, the PM employs the projection matrix to construct $n$-dimensional primitive reciprocal vectors $\{\bm{b}_i^*\}_{i=1}^n$ that span the $n$-dimensional reciprocal space while maintaining integral combination properties. We therefore adopt 
$$\bm{p}_1^*=(1,0),\quad \bm{p}_2^*=(\cos(\pi/6),\sin(\pi/6)),\quad \bm{p}_3^*=(\cos(\pi/3),\sin(\pi/3)),\quad \bm{p}_4^*=(0,1)$$ of $\mathbb{Z}$-rank 4, and the projection matrix 
$$\mathcal{S}=\left(\begin{array}{llll}
1 & \cos (\pi / 6) & \cos (\pi / 3) & 0 \\
0 & \sin (\pi / 6) & \sin (\pi / 3) & 1
\end{array}\right).$$  
As a result, the reciprocal vector for the $n$-dimensional periodic structure is expressed as  
$$
\bm{H} = h_1 \bm{b}_1^* + h_2 \bm{b}_2^* + \cdots + h_n \bm{b}_n^*, \quad h_i \in \mathbb{Z},
$$  
where $\bm{b}_1^*=(1,0,0,0),\bm{b}_2^*=(0,1,0,0),\bm{b}_3^*=(0,0,1,0)$ and $\bm{b}_4^*=(0,0,0,1)$.
Then one can represent any reciprocal vectors $\bm{k} \in \mathbb{R}^d$ of a $d$-dimensional quasicrystal as \cite{chaikin1995principles}
\begin{align*}
\bm{k}=h_1 \bm{b}_1^*+h_2 \bm{b}_2^*+h_3 \bm{b}_3^*+h_4 \bm{b}_4^*, \quad h_i \in \mathbb{Z}.
\end{align*}
Various choices for the coefficient vector 
\[
\bm{h} \triangleq \{h_1, \cdots, h_n\},
\]  
can be employed, 
such as setting some coefficients to zero and imposing constraints on others, leading to different crystal and quasicrystal patterns, as demonstrated in Figure \ref{fig:five_states}.

The Fourier expansion for the $d$-dimensional quasiperiodic function is given by
\begin{equation}
\phi(\bm{r})=\sum_{\bm{H}} \widehat{\phi}(\bm{H}) e^{i\left[(\mathcal{S}\cdot \bm{H})^T \cdot \bm{r}\right]},\quad \bm{r} \in \mathbb{R}^d,\quad \bm{H} \in \mathbb{Z}^n.
\label{eq:PM_solver}
\end{equation}
Denoting by $g_k^T$ the $k$-th row of $\mathcal{S} \cdot \bm{H}$, we have  
$$\mathcal{S} \cdot \bm{H}=\left(\sum_{i=1}^n s_{1 i} \sum_{j=1}^n h_j b_{j i}^*, \cdots, \sum_{i=1}^n s_{d i} \sum_{j=1}^n h_j b_{j i}^*\right)^T \triangleq\left(g_1, \cdots, g_d\right)^T,\quad h_j \in \mathbb{Z},$$
with $b_{j,i}^*, j=1, \cdots, n$ being the $j$-th component of $b_i^*$.
The LP free energy functional then becomes
\begin{equation}
\begin{aligned}
\mathcal{F}(\phi(\bm{r}); c, q, \bmu)= & \frac{1}{2} \sum_{\bm{H}_1+\bm{H}_2=0}\left\{c(1-\sum_{k=1}^d g_k^2)^2(q^2-\sum_{k=1}^d g_k^2)^2-\varepsilon\right\} \widehat{\phi}\left(\bm{H}_1\right) \widehat{\phi}\left(\bm{H}_2\right) \\
& -\frac{\alpha}{3} \sum_{\bm{H}_1+\bm{H}_2+\bm{H}_3=0} \widehat{\phi}\left(\bm{H}_1\right) \widehat{\phi}\left(\bm{H}_2\right) \widehat{\phi}\left(\bm{H}_3\right)\\
&+\frac{1}{4} \sum_{\bm{H}_1+\bm{H}_2+\bm{H}_3+\bm{H}_4=0} \widehat{\phi}\left(\bm{H}_1\right) \widehat{\phi}\left(\bm{H}_2\right) \widehat{\phi}\left(\bm{H}_3\right) \widehat{\phi}\left(\bm{H}_4\right).
\end{aligned}
\label{eq:high_free_energy}
\end{equation}
Substituting Eq. \eqref{eq:PM_solver} into Eq. \eqref{eq:semi_implicit} and using Eq. \eqref{eq:high_free_energy}, one obtains
\begin{equation}
\begin{aligned}
&\left[\frac{1}{\Delta t}+c\left(1-\sum_{k=1}^d g_k^2\right)^2\left(q^2-\sum_{k=1}^d g_k^2\right)^2\right] \widehat{\phi}_{t+\Delta t}(\bm{H}) =\left(\frac{1}{\Delta t}+\varepsilon\right) \widehat{\phi}_t(\mathbf{H})\\
&+\alpha\sum_{\bm{H}_1+\bm{H}_2=\bm{H}} \widehat{\phi}_t\left(\bm{H}_1\right) \widehat{\phi}_t\left(\bm{H}_2\right)-\sum_{\bm{H}_1+\bm{H}_2+\bm{H}_3=\bm{H}} \widehat{\phi}_t\left(\bm{H}_1\right) \widehat{\phi}_t\left(\bm{H}_2\right) \widehat{\phi}_t\left(\bm{H}_3\right),
\end{aligned}
\label{eq:semi_implicit_pm}
\end{equation}
where  $\widehat{\phi}_{t+\Delta t}$ and $\widehat{\phi}_{t}$ represent the Fourier coefficients at time $t+\Delta t$ and $t$, respectively.
A direct evaluation of the convolution terms of \eqref{eq:semi_implicit_pm} are expensive. Instead, one can calculate these nonlinear terms in the physical space and then perform FFT to derive the corresponding Fourier coefficients. Therefore,
the computational complexity of the PM is $$O\left( N_{t}\cdot \calN  \log \calN \right),$$
where $N_{t}$ is the number of time iterations, and $\calN= \left(N_{\bm{H}}\right)^n$ with $N_{\bm{H}}$ being the degrees of freedom of pseudospectral method in each dimension. 

Then we can get the quasicrystal in original space through the projection matrix by \eqref{eq:PM_solver}, and the gradient term $\mathcal{G}(\phi)$ can be computed by
\begin{equation}
\mathcal{G}\left(\phi(\bm{r}\right))=\sum_{\bm{H}}\left[\left(1-\sum_{k=1}^d g_k^2\right)^2\left(q^2-\sum_{k=1}^d g_k^2\right)^2\widehat{\phi}(\bm{H})\right] e^{i\left[(\mathcal{S}\cdot \bm{H})^T \cdot \bm{r}\right]},\quad \bm{r} \in \mathbb{R}^d,\quad \bm{H} \in \mathbb{Z}^n.
\label{eq:PM_solver_grad}
\end{equation}

The field  $\phi(\bm r)$ and corresponding gradient field $\mathcal{G}(\phi(\bm r))$ within any domain can be computed without introducing any Diophantine error by the PM. Following the reference \cite{jiang2014numerical} and \cite{lifshitz1997theoretical}, we choose the edge length of the computational box as  $D=L\times 2\pi$, because the morphologies and the free energy density will not change with denser grids for PM\cite{jiang2014numerical}. For our simulations, we use a time step size of $\Delta t = 0.1$ and a mesh step size of $\Delta x = D/{N_g},$ where $N_g=256$. Then the calculation is performed on $256\times 256$ mesh grids with $L=30$. In this work, we interpret the computational mesh as a unique graph structure,  where the nodes represent the vertices, and the edges define the boundaries of the elements of the triangulation. This conceptual approach is equally applicable to unstructured grids, allowing for flexibility in grid geometry while preserving accuracy.

\section{DiGCA phase classifier}\label{sec:GAPC}

In this section, we first review model order reduction methods in the context of parameterized partial differential equations, motivated by the need to address the computational challenges associated with repeated solving of parameterized LP models for phase diagram construction. 

For the two-dimensional LP model with two distinct length scales, the phase diagram is notably intricate due to the extensive range of parameter values and the multitude of potential stable states. Generating an accurate phase diagram for a broad spectrum of parameters is extremely time-consuming. Ji {\it et al.}\cite{ji2024mcms} introduced the MCMS-RBM, which partitions the parameter domain into multiple components, each designed to simplify and accelerate computations for specific branches of the problem. It significantly accelerates phase diagram generation for LP model, reducing the computation time for producing a detailed phase diagram from several months to just minutes.

\subsection{GCA-ROM for the LP model}
The graph convolutional autoencoder framework for reduced-order modeling (GCA-ROM), introduced in \cite{pichi2024graph}, is a non-intrusive and data-driven approach for nonlinear model order reduction by involving an approximation of the form $u_{\mathcal{N}} \approx \psi(u_N(\bm{\mu}))$ by a nonlinear map $\psi$. The framework employs Graph Neural Networks (GNNs) to encode the reduced manifold, facilitating rapid evaluations of parametrized partial differential equations. 

GNNs are a specialized type of deep learning architecture that excels in extracting valuable information from datasets structured as graphs. They effectively capture various aspects, such as geometric configurations, node relationships, connectivity, and the behavior of features within the graph. Although GNN was originally developed to tackle unstructured grids, we found that it is also effective for handling solutions with distinct structural features. GNNs have emerged as a powerful framework for learning representations of graph-structured data $\mathscr{T} = (\mathcal{V}, \mathcal{E})$, where $\mathcal{V}$ denotes the set of nodes and $\mathcal{E}$ represents edges connecting these nodes. The fundamental operation in GNNs is the \textit{message-passing} mechanism, which propagates information through local neighborhoods of each node $v \in \mathcal{V}$. At the $k$-th layer of a GNN, each node $u \in \mathcal{V}$ aggregates transformed features from its neighborhood $N(u)$ through the following operation:
\begin{equation}
\mathbf{h}_u^{(k)} = \sigma\left(\frac{1}{|N(u)|}\sum_{v \in N(u)} \bm{W}^{(k)}\bm{h}_v^{(k-1)} + \bm{b}^{(k)}\right), \quad k=1,\cdots,K,\quad u\in \mathcal{V}.
\end{equation}
Here, $\bm{h}_u^{(k)} \in \mathbb{R}^{d_k}$ is the hidden state of node $u$ at layer $k$ that is initialized as $\bm{h}_u^{(0)}=u$. The learnable parameters $\bm{W}^{(k)} \in \mathbb{R}^{d_k \times d_{k-1}}$ and $\mathbf{b}^{(k)} \in \mathbb{R}^{d_k}$ are responsible for encoding individualized feature transformations of neighboring nodes. At node $u$, $|N(u)|$ represents the cardinality of its neighborhood. The aggregation across the neighborhood of $u$ facilitates the fusion of information across the graph. Finally, $\sigma(\cdot)$ denotes a typical nonlinear activation function such as ReLU$(\cdot)$, and tanh$(\cdot)$.

A graph convolution network (GCN) \cite{kipf2016semi} is a specific type of GNN aiming at extending the concept of convolutional neural networks (CNNs), operating in regular Euclidean domains, to handle non-grid data. This process is similar to the standard convolutional layers in CNNs, where a fixed filter slides over the pixels of an image to produce gathered information. The primary limitation of CNNs is that their operations lack invariance with respect to the order of nodes. However, this characteristic does not substantially affect their performance when solving PDEs.

GCNs are fundamentally categorized into two principal paradigms: spectral methods\cite{zhu2021simple,li2019specae} and spatial methods\cite{zhou2020fully}. Spectral approaches operate through spectral graph theory by decomposing the graph Laplacian matrix $ L = D - A \in \mathbb{R}^{N \times N},$ where $D$ is the degree matrix and $A$ the adjacency matrix. These methods employ graph Fourier transforms via eigen-decomposition of $L$, establishing rigorous mathematical foundations at the cost of $\mathcal{O}(N^3)$ computational complexity for $N$-node graphs. Spatial methods, in contrast, employ localized neighborhood aggregation mechanisms analogous to conventional CNNs. This paradigm bypasses global graph dependencies through localized sampling operations, achieving $\mathcal{O}(|\mathcal{E}|)$ computational complexity relative to edge count $|\mathcal{E}|$. 

The MoNet framework \cite{monti2017geometric} is a spatial method that performs geometrically informed convolutions on non-Euclidean domains using Gaussian mixture models. To this end, it introduces edge-specific pseudo-coordinates $e_{uv} \in \mathbb{R}^d$ , typically Euclidean distances between nodes $u$ and $v$, to parameterize Gaussian kernels:

\begin{equation}
    \omega^q(e_{uv}) = \exp\left(-\frac{1}{2}(e_{uv} - \mu_q)^\top \Sigma_q^{-1}(e_{uv} - \mu_q)\right)
\end{equation}
where $\{\mu_q \in \mathbb{R}^d, \Sigma_q \in \mathbb{R}^{d \times d}\}_{q=1}^Q$ are trainable parameters for $Q$ filters. The node update rule combines spatial and feature information through:

\begin{equation}
\label{eq:MoNet}
    \bm{h}_u^{(k+1)} = \sigma\left(\frac{1}{|N(u)|}\sum_{v \in \mathcal{N}(u)} \sum_{q=1}^Q {\bf \omega}^q(e_{uv}) \odot \bm{W}_q \bm{h}_v^{(l)}\right)
\end{equation}
where $\bm{ W}_q \in \mathbb{R}^{d \times d}$ denotes learnable weight matrices. This formulation introduces geometric bias while maintaining linear complexity relative to edge count $|\mathcal{E}|$.

The GCA-ROM framework operates through an autoencoder-decoder architecture that implements a convolve-then-reshape paradigm, performing convolution directly on the original geometric or grid structure. Autoencoders are unsupervised learning models that learn low-dimensional representations of high-dimensional data. By mapping the high-dimensional PDE solutions into a low-dimensional latent space, autoencoders can efficiently construct reduced-order models. Subsequently, a secondary network is trained to predict the evolution of the solution within this latent space, allowing the online computations to be performed entirely in the low-dimensional space. This preserves the data's inherent geometric relationships and spatial dependencies, effectively capturing the inductive bias or structural prior embedded in the data. In contrast, standard CNNs use a reshape-then-convolve strategy, which reshapes the data before applying convolution. This pre-processing step can disrupt the original structure, leading to a loss of critical geometric and spatial information, making GCA-ROM better suited for tasks involving irregular geometries or graph-structured data.

Let $\mathcal{G}=(\mathcal{V},\mathcal{E})$ define the graph structure with nodal features $\bm{h}^{(0)} \in \mathbb{R}^{|\mathcal{V}| \times d}$ and parameterized dataset $\Xi = \{(\bm{\mu}^i, u_{\mathcal{N}}(\bm{\mu}^i))\}_{i=1}^{N_s}$, where $u_{\mathcal{N}}(\bm{\mu}^i) \in \mathbb{R}^{|\mathcal{V}| \times d}$ denotes high-fidelity solutions of parametrized PDEs. The autoencoder $\mathcal{I}_W: \Omega_{\mathcal{N}} \to u_{\mathcal{N}}$ implements dimension reduction through symmetric geometric operations:
\begin{equation}
    \mathcal{I}_W = \psi_W \circ \phi_W \quad \text{where} \quad 
    \begin{cases}
    \text{Encoder } \phi_W: & \quad \mathbb{R}^{|\mathcal{V}| \times d} \xrightarrow{\mathcal{GC}} \mathbb{R}^{|\mathcal{V}| \times r} \xrightarrow{\mathcal{M}} \mathbb{R}^N ,\\
    \text{Decoder } \psi_W: & \quad \mathbb{R}^N \xrightarrow{\mathcal{M}^\dagger} \mathbb{R}^{|\mathcal{V}| \times r} \xrightarrow{\mathcal{GC}^\dagger} \mathbb{R}^{|\mathcal{V}| \times d}.
    \end{cases}
\end{equation}
The MoNet geometric convolution operator $\mathcal{GC}$ \cite{monti2017geometric}, defined in Eq.~\eqref{eq:MoNet}, performs a localized feature transformation while preserving nodal and edge attributes through edge-conditioned weight modulation. Subsequent dimension reduction is achieved by the bottleneck map $\mathcal{M}: \mathbb{R}^{|\mathcal{V}| \times r} \rightarrow \mathbb{R}^N$ with $\text{rank}(\mathcal{M}) = N \ll |\mathcal{V}|$, which learns the latent behavior in a vector. The decoder phase utilizes the pseudoinverse projection $\mathcal{M}^\dagger: \mathbb{R}^N \rightarrow \mathbb{R}^{|\mathcal{V}| \times r}$ to restore nodal dimensions via convex combination of latent features, paired with transposed convolution $\mathcal{GC}^\dagger$ that reverses the feature transformation through parameter-sharing transpose operations. This operator quartet maintains geometric consistency through $\mathcal{GC}$-$\mathcal{GC}^\dagger$ duality and $\mathcal{M}$-$\mathcal{M}^\dagger$ pseudoinverse relationships, ensuring topological preservation across scale transitions.

The encoder-decoder structure is designed to approximate the identity operator by leveraging manifold learning. Specifically, the composite operator satisfies the following:  
\begin{equation}  
    \mathcal{I}_W = \psi_W \circ \phi_W \approx \mathcal{I}_{|\mathcal{V}| \times d}.  
\end{equation}  
The unsupervised learning objective enforces the following minimization problem:
\begin{equation}  
    \min_W \mathcal{L}_s(W), \quad {\rm with } \quad \mathcal{L}_s(W) \triangleq \| \psi_W(\phi_W(u)) - u \|_{2}^2,  
\end{equation}  
Here, $\mathcal{L}_s(W)$ measures the solution reconstruction loss.

 By combining the structural capabilities of GNNs with the dimensionality reduction power of autoencoders, GCA-ROM effectively addresses the challenges posed by complex geometries, irregular grids, and high-dimensional solution spaces. Upon encoding the latent variables with the autoencoder, we utilize this dataset to carry out the supervised learning task using the multi-layer perceptron (MLP). In doing so, we non-intrusively establish the mapping $\hat{u}_{\mathcal{N}}(\bm{\mu})=\text{MLP}(\bm{\mu})$, thereby enabling the recovery of the reduced coefficient. The latent variable representations from both paths are constrained through minimizing the latent variable reconstruction loss:
$$
\min_W \mathcal{L}_v(W), \quad {\rm with } \quad \mathcal{L}_v(W) \triangleq \| \phi_W(u_{\mathcal{N}}(\bm{\mu})) - \text{MLP}(\bm{\mu}) \|_{2}^2,
$$
where $\phi_W(u_{\mathcal{N}}(\bm{\mu}))$ denotes the encoder output.

The two minimization problems can be weighted by the hyperparameter $\lambda$, thus constructing the final loss function:
\begin{equation}
\label{eq:gca-loss}
\mathcal{L} = \underbrace{\| \psi_W(\phi_W(u_{\mathcal{N}})) - u_{\mathcal{N}} \|_{2}^2}_{\text{Solution reconstruction}} + \lambda \underbrace{ \| \phi_W(u_{\mathcal{N}}(\bm{\mu})) - \text{MLP}(\bm{\mu}) \|_{2}^2}_{\text{Latent variable reconstruction}}.
\end{equation}

\begin{algorithm}
\caption{DiGCA Offline Training Procedure}
\label{alg:gcaron}

\textbf{Input:} 
Training dataset $\Xi_{\text{train}}=\{u_{\mathcal{N}(\bm{\mu}^i)},\Omega_{\mathcal{N}(\bm{\mu}^i)}\}_{i=1}^{N_s}$ 

\textbf{Initialize:}
Randomly initialize encoder $\phi_W$, decoder $\psi_W$ and bottleneck map MLP

\textbf{For }$k=1$ \textbf{to} $K_{\text{max}}$:
\begin{algorithmic}[1]
\State Compute latent representation: $u_N(\bm{\mu}^i) = \phi_W(u_{\mathcal{N}(\bm{\mu}^i)}), \quad i \in \Xi_{\text{train}}$
\State Multi-layer perceptron: $\widehat{u_N}(\bm{\mu^i}) = \text{MLP}(\bm{\mu^i}), \quad i \in \Xi_{\text{train}}$
\State Reconstruct graph data: $\widehat{u_{\mathcal{N}}(\bm{\mu}^i)} = \psi_W(u_N(\bm{\mu}^i) ), \quad i \in \Xi_{\text{train}}$
\State Evaluate loss function:
$$
\mathcal{L} = \frac{1}{N_s}\sum_{i=1}^{N_s} \left(\underbrace{\| \psi_W(\phi_W(u_{\mathcal{N}}(\bm{\mu}^i))) - u_{\mathcal{N}} (\bm{\mu}^i)\|_{2}^2}_{\text{Solution reconstruction}} + \lambda \underbrace{ \| \phi_W(u_{\mathcal{N}}(\bm{\mu}^i)) - \text{MLP}(\bm{\mu}^i) \|_{2}^2}_{\text{Latent variable reconstruction}} \right)
$$
\State Update the parameters of $\phi_W$, $\psi_W$, and $\text{MLP}$ through backpropagation.
\end{algorithmic}

\end{algorithm}

\subsection{Training loss with high-order derivatives}
\label{sec:gca-rom}

However, accurately computing derivatives has long been a challenge in using neural networks to solve PDEs\cite{lyu2022mim}. Higher-order derivatives, which require multiple differentiations, often result in significant errors, particularly when calculating energy. To address this issue, we incorporate derivative information into the network's training loss function. We consider the graph dataset $\Xi_{\text{train}}=\{u_{\mathcal{N}(\bm{\mu}^i)},\mathcal{G}(u_{\mathcal{N}(\bm{\mu}^i)}),\Omega_{\mathcal{N}(\bm{\mu}^i)}\}_{i=1}^{N_s}$ which consists of $N_s$ solutions $u_{\mathcal{N}(\bm{\mu}^i)}$ and nonlocal term 
 $\mathcal{G}(u_{\mathcal{N}(\bm{\mu}^i)})$ of a parameterized PDE defined over computation domain, corresponding to the parameter set $\{\mu^i\}_{i=1}^{N_s}$. For simply, we define $\mathrm{U}_{\mathcal{N}(\bm{\mu}^i)} = [u_{\mathcal{N}(\bm{\mu}^i)}, 
 \lambda_u\mathcal{G}(u_{\mathcal{N}(\bm{\mu}^i)})]$, where $\lambda$ is a hyperparameter can balance the loss between solutions and gradients. 

The architecture for offline training consists of an autoencoder to encode the information into a low-dimensional space and a non-intrusive MLP to map the parameter with the latent vector $\mathrm{U}_N(\bm{\mu}) = \text{MLP}(\bm{\mu})$. The latent vector also obtained by encoding the solution field $\hat{\mathrm{U}}_N(\bm{\mu})=\Phi(\mathrm{U}_{\mathcal{N}}(\bm{\mu}))$. So we define the first term in the loss function for the learning task as
\begin{equation}
\mathcal{L}_v = \frac{1}{N_{\text{train}}}\sum_{i=1}^{N_{\text{train}}}\lVert \mathrm{U}_N(\bm{\mu}^i)- \hat{\mathrm{U}}_N(\bm{\mu}^i) \rVert _2^2.
\end{equation}
Then the decoding structure mirrors the encoding process, using the same operations but in reverse order. So we can reconstruct the solution by $\hat{U}_{\mathcal{N}}(\bm{\mu})=\Psi(\hat{U}_N(\bm{\mu}))$ and defined the second term of the loss function as
\begin{equation}
\mathcal{L}_s = \frac{1}{N_{\text{train}}}\sum_{i=1}^{N_{\text{train}}}\lVert \mathrm{U}_{\mathcal{N}}(\bm{\mu}^i)- \hat{\mathrm{U}}_{\mathcal{N}}(\bm{\mu}^i) \rVert _2^2.
\end{equation}
The total loss functions can be balanced through a hyperparameter $\lambda$ as follows:
\begin{equation}
\mathcal{L} = \lambda \mathcal{L}_{v} +\mathcal{L}_s. 
\end{equation}

In the online phase, we only use the trained MLP to evaluate the $\mathrm{U}_N(\bm{\mu}) = \text{MLP}(\bm{\mu})$ for a new parameter and then decompress through the graph decoder to recover the solution over its geometry. The method can directly compute the components $E_1$ and $E_2$ for all parameters and all five states based on Eq.~\eqref{eq:GH-defined}. According to free energy theory, the solution with the minimum energy corresponds to the equilibrium state, which is defined as the phase associated with specific material parameters in the system.

\begin{figure}[htbp]
    \centering    
    \includegraphics[width=1\textwidth]{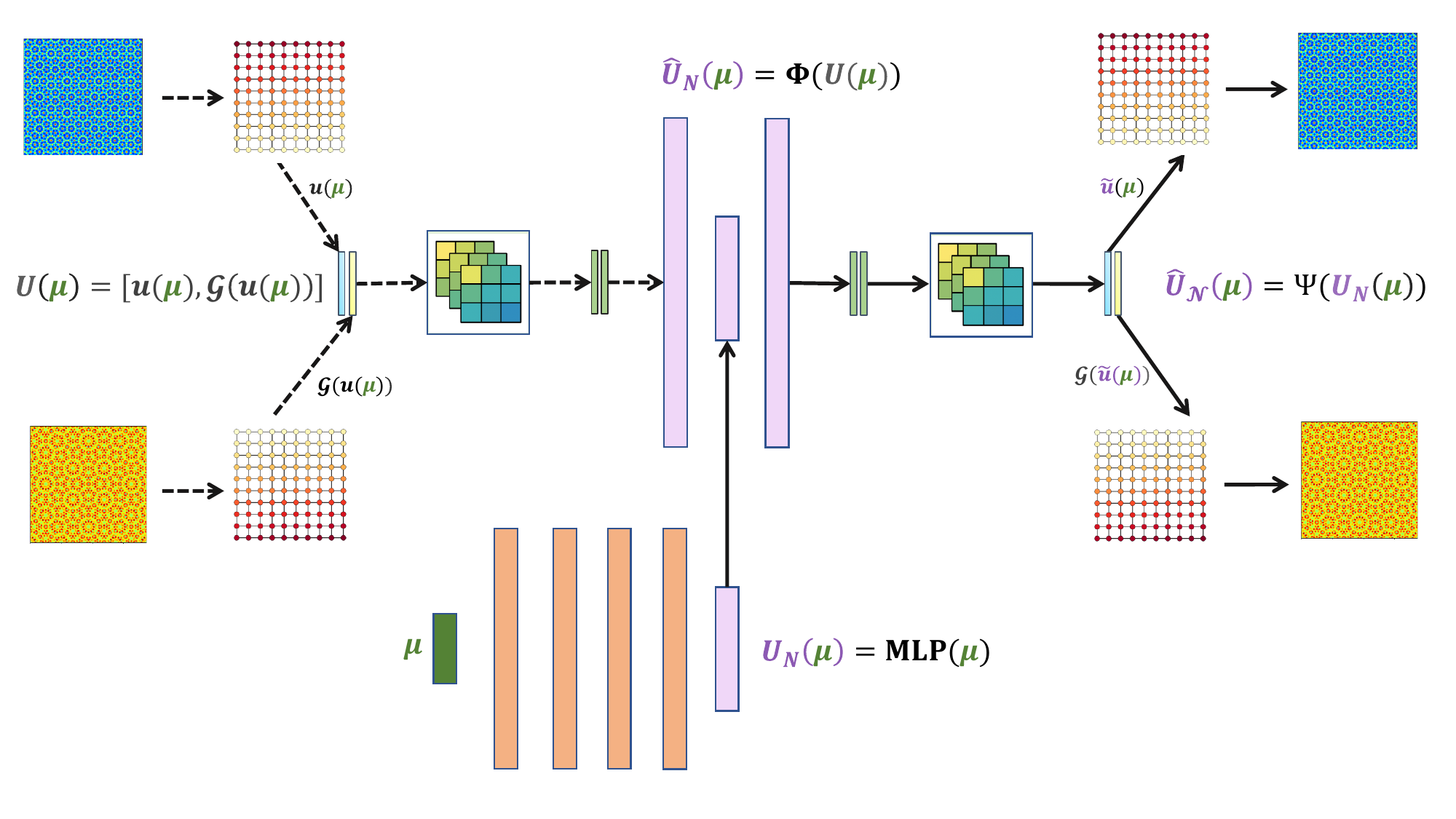}
    \caption{DiGCA schematic for Lifshitz-Petrich model.}    \label{fig:ML_LP_offline1}
\end{figure}

The architecture of DiGCA is shown \Cref{fig:ML_LP_offline1}.

\begin{algorithm}
\caption{DiGCA Offline Training Procedure}
\label{alg:gcaron}

\textbf{Input:} 
Training dataset $\Xi_{\text{train}}=\{u_{\mathcal{N}(\bm{\mu}^i)},\mathcal{G}(u_{\mathcal{N}(\bm{\mu}^i)}),\Omega_{\mathcal{N}(\bm{\mu}^i)}\}_{i=1}^{N_s}$ 

\textbf{Initialize:}
Randomly initialize encoder $\Phi$, decoder $\Psi$ and bottleneck map MLP

\textbf{For }$k=1$ \textbf{to} $K_{\text{max}}$:
\begin{algorithmic}[1]
\State Compute latent representation: $U_N(\bm{\mu}^i) = \Phi(U_{\mathcal{N}(\bm{\mu}^i)}), \quad i \in \Xi_{\text{train}}$
\State Multi-layer perceptron: $\widehat{U_N}(\bm{\mu^i}) = \text{MLP}(\bm{\mu^i}), \quad i \in \Xi_{\text{train}}$
\State Reconstruct graph data: $\widehat{U_{\mathcal{N}}(\bm{\mu}^i)} = \Psi_W(U_N(\bm{\mu}^i) ), \quad i \in \Xi_{\text{train}}$
\State Evaluate loss function:
$$
\mathcal{L} = \frac{1}{N_s}\sum_{i=1}^{N_s} \left(\underbrace{\| \Psi_W(\Phi_W(U_{\mathcal{N}}(\bm{\mu}^i))) - U_{\mathcal{N}} (\bm{\mu}^i)\|_{2}^2}_{\text{Solution reconstruction}} + \lambda \underbrace{ \| \Phi_W(U_{\mathcal{N}}(\bm{\mu}^i)) - \text{MLP}(\bm{\mu}^i) \|_{2}^2}_{\text{Latent variable reconstruction}} \right)
$$
\State Update the parameters of $\phi_W$, $\Psi_W$, and $\text{MLP}$ through backpropagation.
\end{algorithmic}
\end{algorithm}

We employ two tricks in the network architecture and the training process. The first trick is the activation function. We use $\sin(w \cdot x)$ as the activation function, where $w$ is a trainable parameter. This activation function effectively captures high-frequency oscillations induced by higher-order derivatives, making it particularly suitable for problems where the solution exhibits rapid variations. Numerical results demonstrate that this choice confirms that this choice outperforms conventional activation functions such as $\text{Tanh}(\cdot)$ and $\text{Relu}(\cdot)$ in terms of capturing these features.  The second trick is minibatch training. The minibatch training processes small subsets of the dataset, known as minibatches, during each training iteration. Unlike batch training, which uses the entire dataset at once, or stochastic training, which processes one data point at a time, minibatch training offers a practical middle ground. It updates model parameters more frequently than batch gradient descent, enabling faster progress, while being more stable than the noisy updates of stochastic gradient descent. This balanced approach not only accelerates convergence but also maintains a good level of stability, making it an effective and widely used training strategy.

\subsection{Multi-state phase classifier}
\label{sec:classifier}
During the traditional procedure for any given parameter value $\mu$, the LP model should be solved several times with different initial value given in the corresponding candidate state. So to resolve the phase diagram even on a relatively small parameter domain, thousands of simulations and free energy computation are needed. 
The structures and properties of solutions in different states exhibit significant differences, so Ji {\it et al} \cite{ji2024mcms} developed the multi-component multi-state reduced basis method, named MCMS-RBM. The generic framework of MCMS-RBM consists of various components, each serving to reduce a specific problem branch associated with a particular part of the parameter domain. Drawing inspiration from this concept, our objective is to develop a pre-trained subnet with multi-components that can efficiently predict solutions for parameterized problems by DiGCA and design a phase classifier by a deep neural network, which can predict the phase for any parameter precisely and quickly.

\begin{figure}[htbp]
    \centering
    \includegraphics[width=\textwidth]{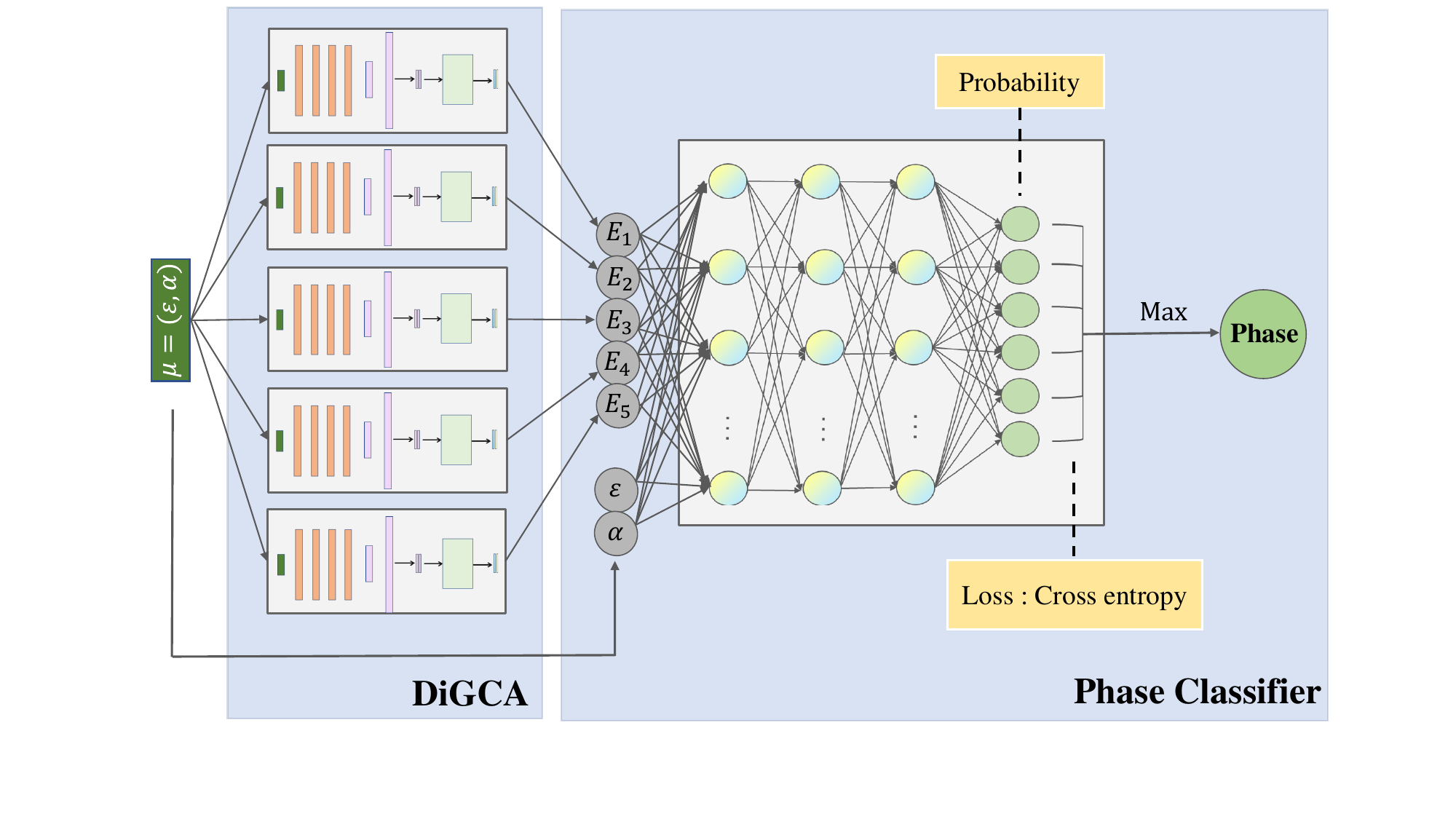}
    \caption{The schematic of DiGCA with phase classification.}
    \label{fig:ML_LP_online}
\end{figure}

The schematic of our design is provided in Figure \ref{fig:ML_LP_online}. Every component of the subnet is connected with the stable order parameter in LP model, which has the minimum energy in Eq.~\eqref{eq:phi_energy}. Thus we should divide all the parameters and solutions in the dataset into five parts by the stable state and train the network separately. If we have a set of stable solutions corresponding to specific parameters, cluster analysis methods such as k-means can be used to classify the dataset. Additionally, images or data obtained from experiments can be processed and incorporated into the dataset as part of the training set. In fact, the network can even be trained using only experimentally obtained data, which would significantly reduce the offline cost of generating high-accuracy solutions.

Using the five subnets of the DiGCAs network, we obtain energy features corresponding to a new parameter $\bm{\mu} = (\varepsilon,\alpha)$. Different from minimizing the free energy to find the stable state, the neural network phase classifier for LP model uses $7$ features 
\[\mathcal{E}=\{\varepsilon,\alpha,\mathrm{E}^{\mathrm{ QC}},\mathrm{E}^{\mathrm{C6}},\mathrm{E}^{\mathrm{LQ}},\mathrm{E}^{\mathrm{T6}},\mathrm{E}^{\mathrm{Lam}}\},\]
from the five DiGCA components and the labels are the stable state corresponding to $\bm{\mu}$. The 7-dimensional feature vector $\mathcal{E} \in \mathbb{R}^7$ effectively captures the system's energy characteristics while achieving substantial memory reduction compared to the full-order parameterization with $\mathcal{N}$ degrees of freedom. Each phase label $\text{S} \in \{\text{QC}, \text{C6}, \text{LQ}, \text{T6}, \text{Lam}, \text{Lq}\}$ is encoded into a canonical basis vector $\mathcal{P}_{\text{S}} \in \mathbb{R}^6$ via one-shot encoding:

\begin{equation}
    \mathcal{P}_{\text{S}} = [\delta_{\text{QC},\text{S}}, \delta_{\text{C6},\text{S}}, \delta_{\text{LQ},\text{S}}, \delta_{\text{T6},\text{S}}, \delta_{\text{Lam},\text{S}}, \delta_{\text{Lq},\text{S}}]^\top,
\end{equation}
where $\delta_{i,j}$ denotes the Kronecker delta function. This encoding scheme guarantees dimensional consistency with the neural network's 6-dimensional output layer. The training dataset $\Xi_p = \{(\mathcal{T}_i, \mathcal{P}_i)\}_{i=1}^N$ comprises input vectors $\mathcal{T}_i \in \mathbb{R}^7$ containing physical parameters paired with corresponding label vectors $\mathcal{P}_i \in \mathbb{R}^6$, where each label vector represents the phase probability distribution in six-dimensional space.

The phase classifier is a fully connected neural network with layer widths [7,40,40,40,6] and $\text{Tanh}(\cdot)$ as the activation function. The trained classifier implements a probabilistic mapping:
\begin{equation}
 \bm{\mu} \mapsto [P_{\text{QC}}(\bm{\mu}), P_{\text{C6}}(\bm{\mu}), P_{\text{LQ}}(\bm{\mu}), P_{\text{T6}}(\bm{\mu}), P_{\text{Lam}}(\bm{\mu}), P_{\text{Lq}}(\bm{\mu})] \in [0,1]^6,
\end{equation}
with the phase-determination criterion:
\begin{equation}
    \text{Phase}(\bm{\mu}) = \mathop{\arg\max}\limits_{\text{S} \in \{\text{QC}, \text{C6}, \text{LQ}, \text{T6}, \text{Lam}, \text{Lq}\}} P_{\text{S}}(\bm{\mu}).
\end{equation}
The loss function used is the cross-entropy loss between the labels and the predictions. During the offline phase, we train the neural network using the Adam optimizer for $3000$ epochs. In the online stage, the trained classifier can evaluate the parameter phase instantaneously.

\section{Numerical results}

\label{sec:numerics}

In this section, we test the proposed DiGCA on the two-dimensional quasiperiodic LP model parameterized by the reduced temperature $\varepsilon$ and the phenomenological parameter $\alpha$ delineating the level of asymmetry. Furthermore, we highlight its efficiency and accuracy by adopting the adaptive phase diagram generation algorithm to produce a phase diagram that is as accurate as what's produced by a full order model.

The parameter domain is set to be $\mathcal{D}_{\Xi}=[-0.01,0.05]\times [0,1]$. Due to the multi-phase nature of the underlying physical problem, $\mathcal{D}_{\Xi}$ is inherently composed of several branches, each corresponding to a distinct stable state of the system. For the LP model, these branches naturally partition into six categories: QC, C6, LQ, T6, Lam and liquid (Liq). For each branch, we randomly select $N_s=200$ parameter samples to construct the dataset. This dataset is then split into a training set and a validation set according to a ratio of $r_t = 75\%$, denoted by 
$\mathcal{D}_{\Xi}^{\mathrm{train}}$ and $\mathcal{D}_{\Xi}^{\mathrm{test}}$, respectively. Specifically, 
$r_t N_s$ samples are used to train the corresponding subnet, while the remaining samples $(1-r_t) N_s$ serve as the validation set. Although the phase diagram can be constructed using 1000 training samples, the resulting phase boundaries remain considerably coarse-particularly in regions near the intersections of the QC, LQ, and C6 phases. This observation highlights the necessity for efficient algorithms capable of real-time phase diagram prediction. Such methods should offer high-resolution boundary delineation while significantly reducing computational cost, thereby making them suitable for practical applications that demand rapid evaluation.

We set the degree of freedom of the Fourier spectral method in PM in each direction as $N_{\mathbf{H}}=32$ to solve the LP model. The complexity for every parameter is $O(N_t\cdot N_{\mathbf{H}}^4\log{N_{\mathbf{H}}})$, where $N_t$ is the number of time iterations. The dataset for DiGCA also needs to discretize the parameter set $\mathcal{D}_{\Xi}.$ In the LP model, we use uniform grid with a mesh step size of $\Delta x = D/{N_g}$ to generate the high fidelity solutions, where $N_g=256$.

\begin{figure}[htbp]
    \centering
\includegraphics[width=0.48\linewidth]{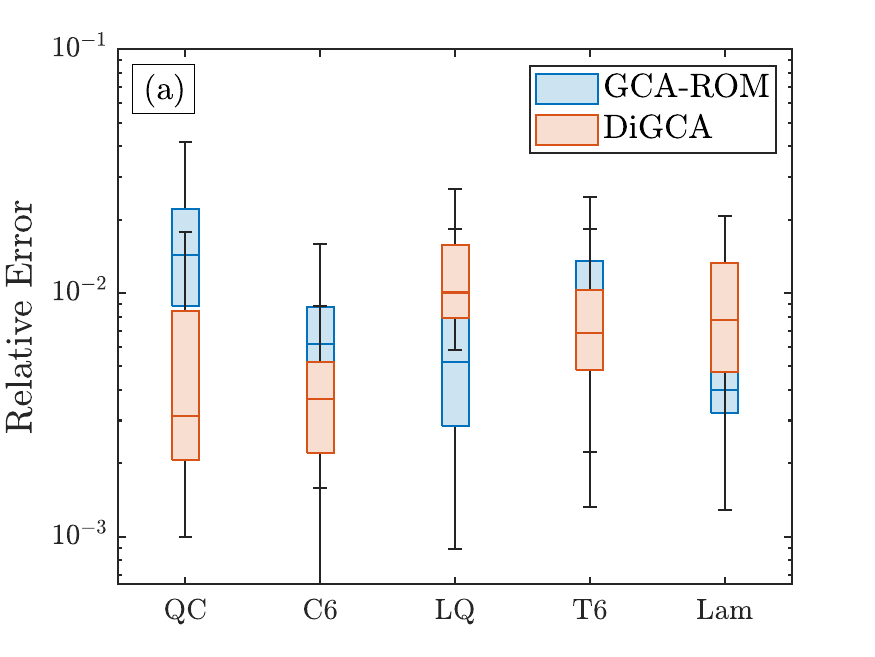}
\includegraphics[width=0.48\linewidth]{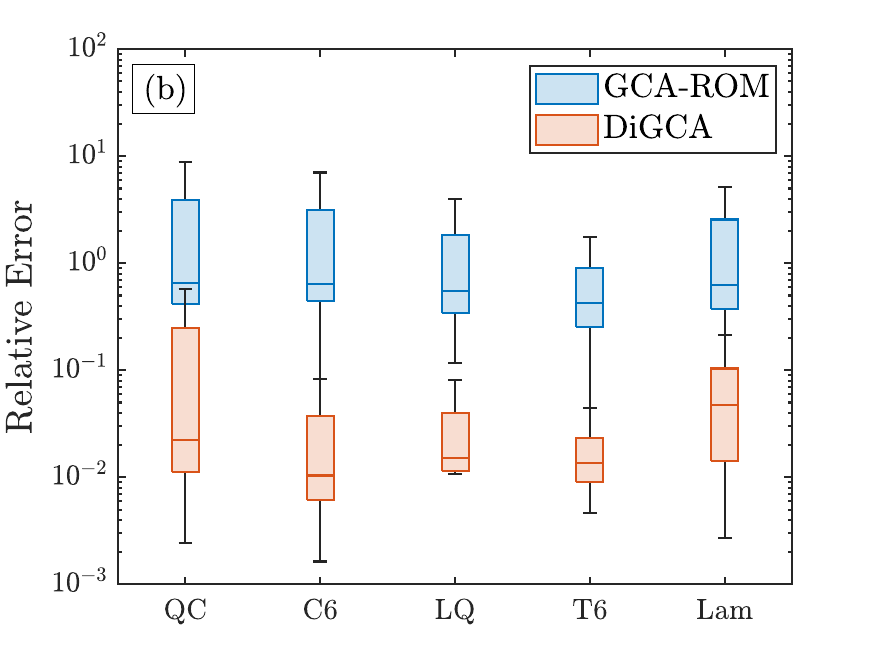}
\caption{Comparison of the relative errors between the GCA-ROM and DiGCA: (a) the solution and (b) the nonlocal term.}
\label{fig:ROM_err}
\end{figure}

We compare two methodologies: the regular graph convolutional autoencoder, which trains the model using only the solution, and our proposed derivative-informed graph convolutional autoencoder, which integrates both the solution and its gradient during training. As shown in \Cref{fig:ROM_err}, we assessed the neural network's performance on $\mathcal{D}_{\Xi}$ for each component by examining relative $L_2$ errors. The derivative-informed training method significantly improves gradient evaluation accuracy while achieving solution prediction performance comparable to the regular method.

\begin{figure}
    \centering    \includegraphics[width=1.05\linewidth]{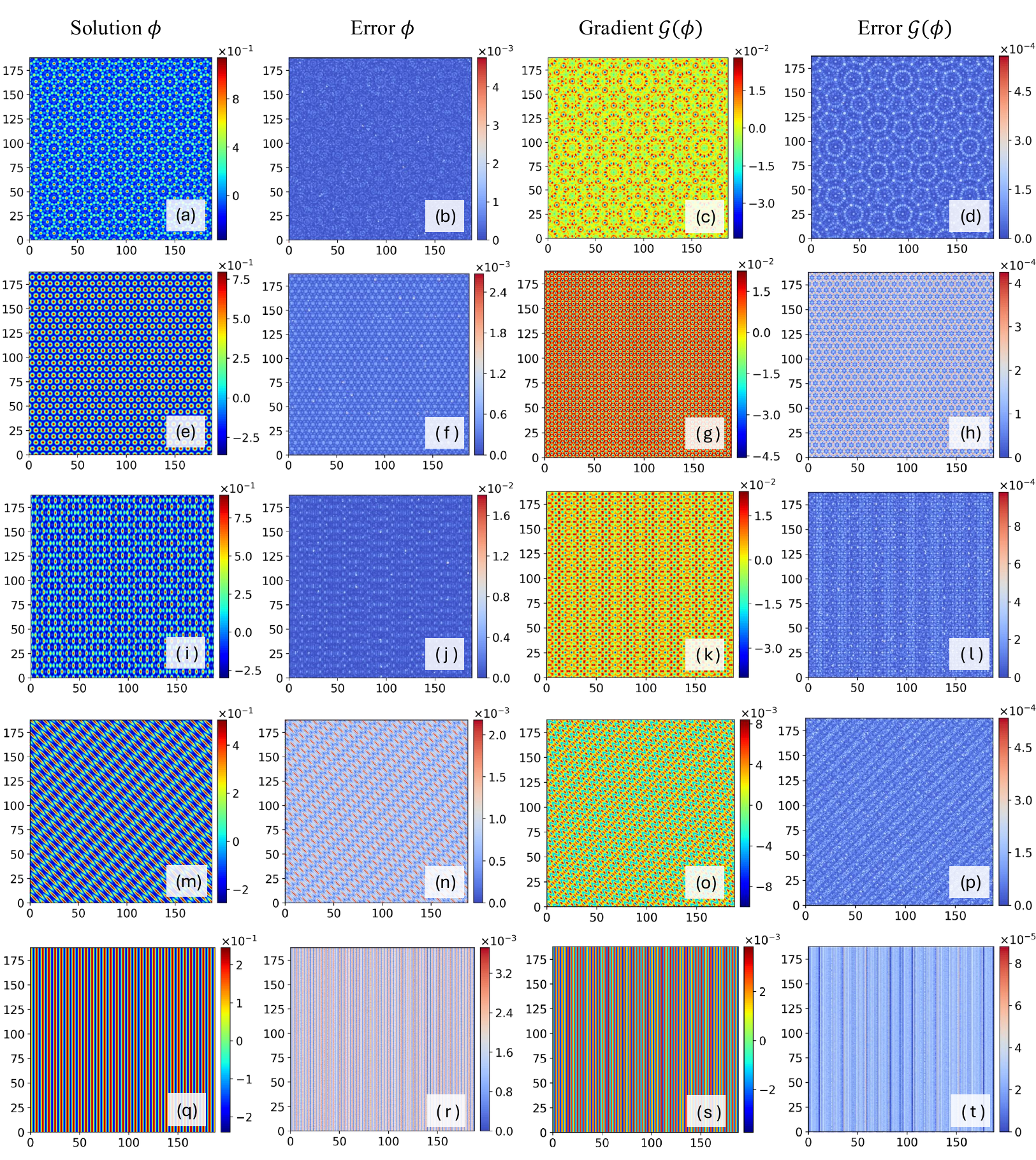}
    \caption{Derivative-informed graph convolutional autoencoder results. Columns (left to right): Order parameter solution ($\phi$), solution error, gradient term ($\mathcal{G}(\phi)$), and gradient term error. Rows correspond to different parameter states: (a-d) QC: $\bmu=(0.0011,0.8300)$, (e-h) C6: $\bmu=(0.0396,0.7740)$, (i-l) LQ: $\bmu=(0.0083,0.8020)$, (m-p) T6 :$\bmu=(0.0477,0.3560)$, and (q-t) Lam: $\bmu=(0.0431, 0.1080)$. }
    \label{fig:GCA-solution}
\end{figure}

Then we present the solution $\phi(\bm{r})$ and nonlocal term $\mathcal{G}(\phi(\bm{r}))$ computed by the enhanced graph convolutional autoencoder on the testing set in \Cref{fig:GCA-solution}. The pointwise relative errors for both the solution and nonlocal term remain mostly below 5\%. 
The precision attained in simultaneously capturing both the primary solution and its nonlocal counterpart underscores the robustness of our proposed methodology, suggesting its potential for broader applications in nonlinear problems involving nonlocal operators. The sub-percent level errors observed throughout the computational domain validate the autoencoder's ability to maintain physical consistency between the solution and its derived quantities.

\begin{figure}[htbp]
\centering
\includegraphics[width=0.48\linewidth]{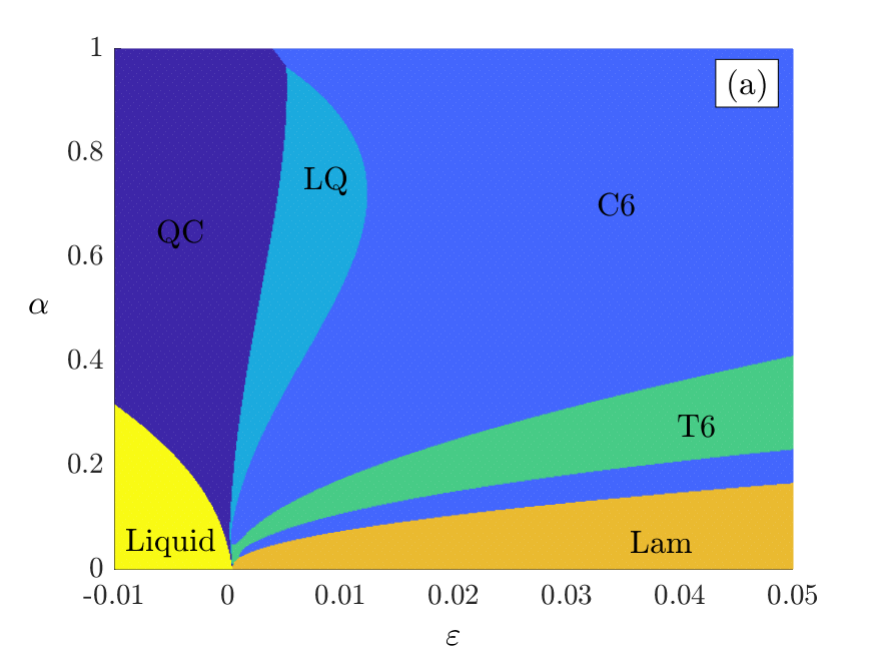}
\includegraphics[width=0.48\linewidth]{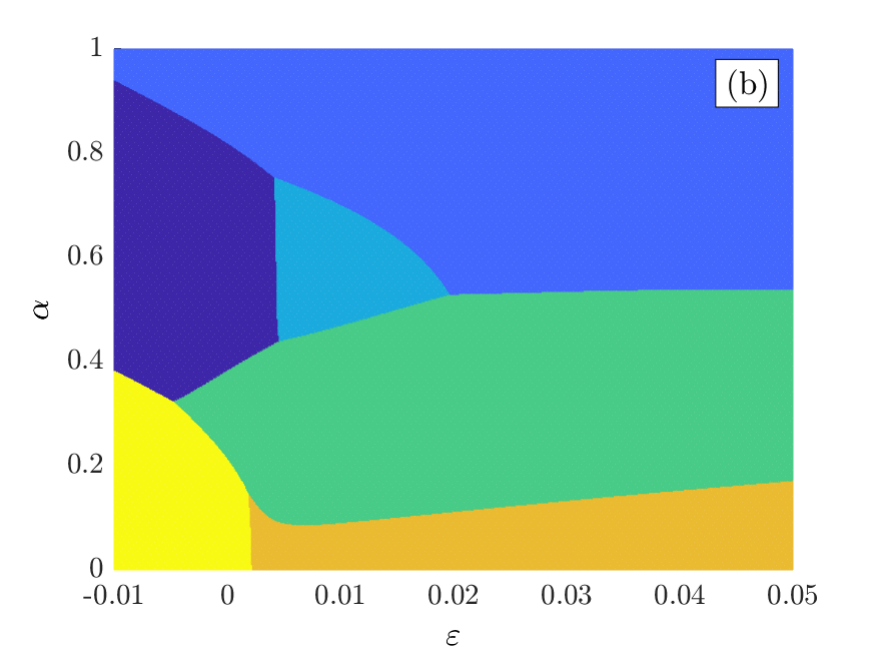}
\includegraphics[width=0.48\linewidth]{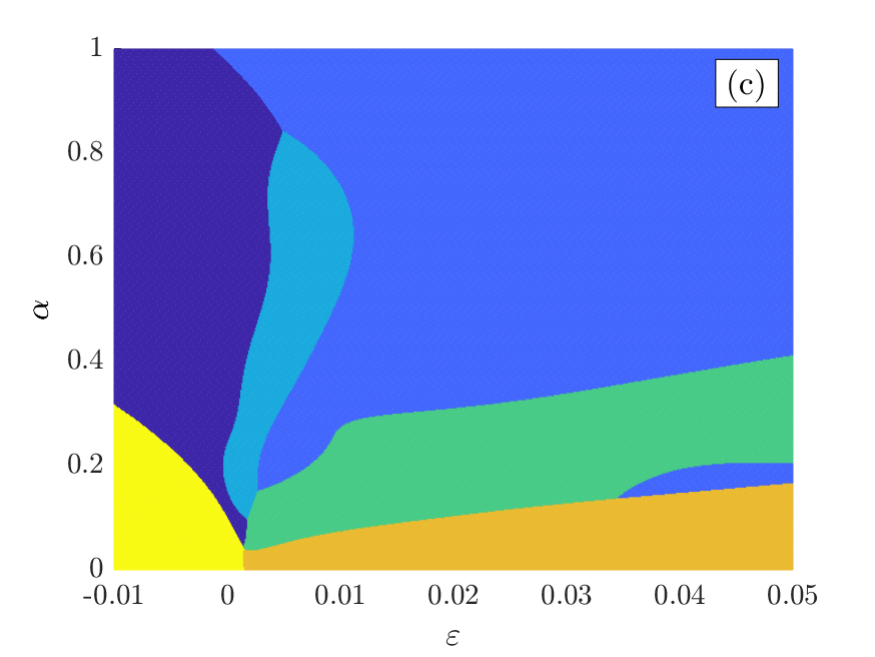}
\includegraphics[width=0.48\linewidth]{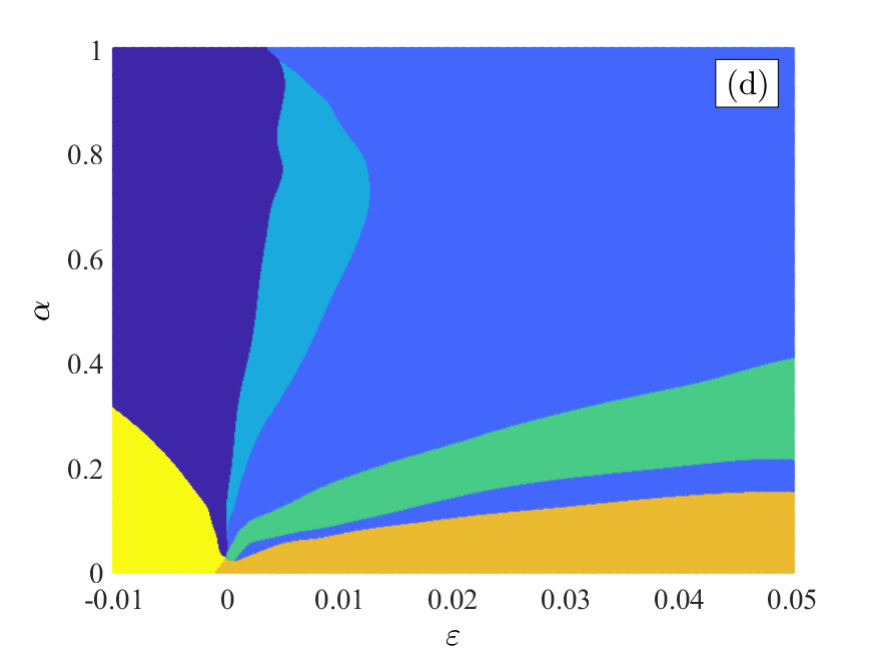}
\caption{Phase diagrams generated by different methods. (a): MCMS-RBM, (b): GCA-ROM, (c): DiGCA, (d): DiGCA with phase classifier.}
\label{fig:ROM_phase}
\end{figure}

Furthermore, our comparative analysis reveals striking differences in phase boundary detection capability between conventional and the derivative-informed approach. As shown in Figure \ref{fig:ROM_phase} (a-c), the conventional method (b) produces results that significantly deviate from the reference solution (a). It reveals rough classification trends rather than well-defined boundaries. In stark contrast, our DiGCA (c) successfully reconstructs most phase boundaries with high accuracy, including their characteristic curvatures and topological connections. While some minor discrepancies remain near the origin where all stable states meet, the improved method demonstrates superior capability in resolving fine boundary features compared to the conventional approach. The results clearly demonstrate that gradient-enhanced learning enables reliable phase diagram prediction where conventional methods do not capture detailed boundary information.

To enhance the precision of phase boundary determination, we implement a deep neural network classifier to post-process the solutions generated by our graph convolutional autoencoder. The classification network, with architecture [7,40,40,40,6], is trained offline using the complete feature-label pairs from dataset $\Xi$. This two-stage approach combines the representational power of enhanced graph convolutional autoencoder for solution generation with the discriminative capability of deep neural networks for precise phase classification, as shown in \Cref{fig:ROM_phase} (d). The 7-dimensional input features capture the essential physical parameters, while the 6 output nodes correspond to distinct phases in the system. Through this hierarchical learning framework, we achieve improved resolution of phase boundaries compared to the former methods, particularly in regions where traditional approaches might show inconsistency.

\begin{figure}
    \centering
    \includegraphics[width=0.32\linewidth]{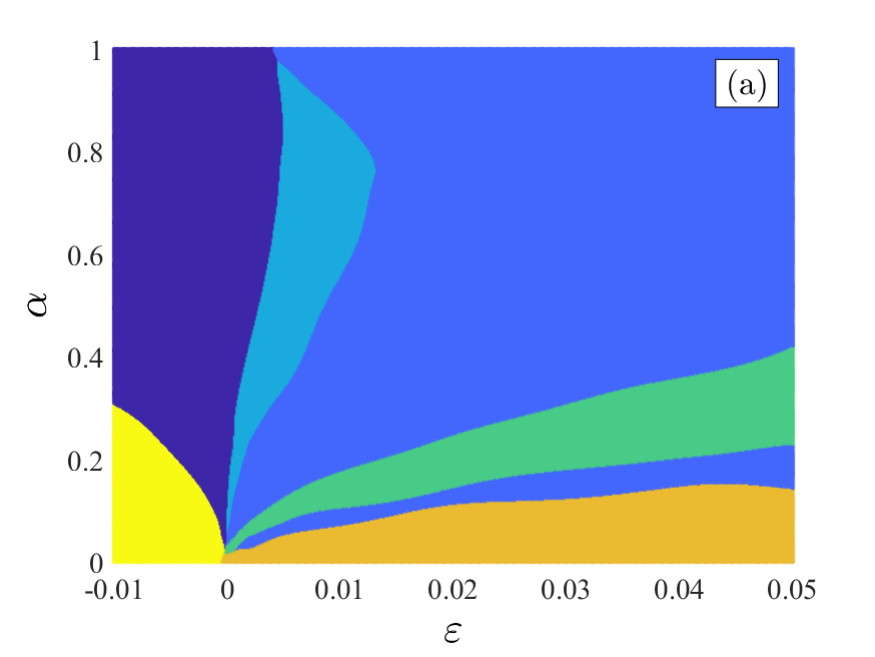}
    \includegraphics[width=0.32\linewidth]{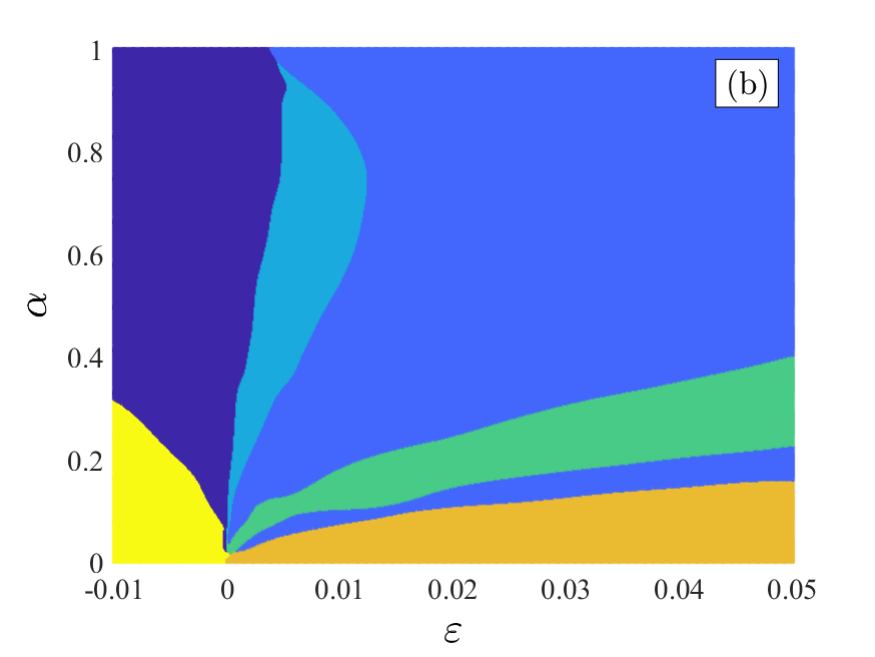}
    \includegraphics[width=0.32\linewidth]{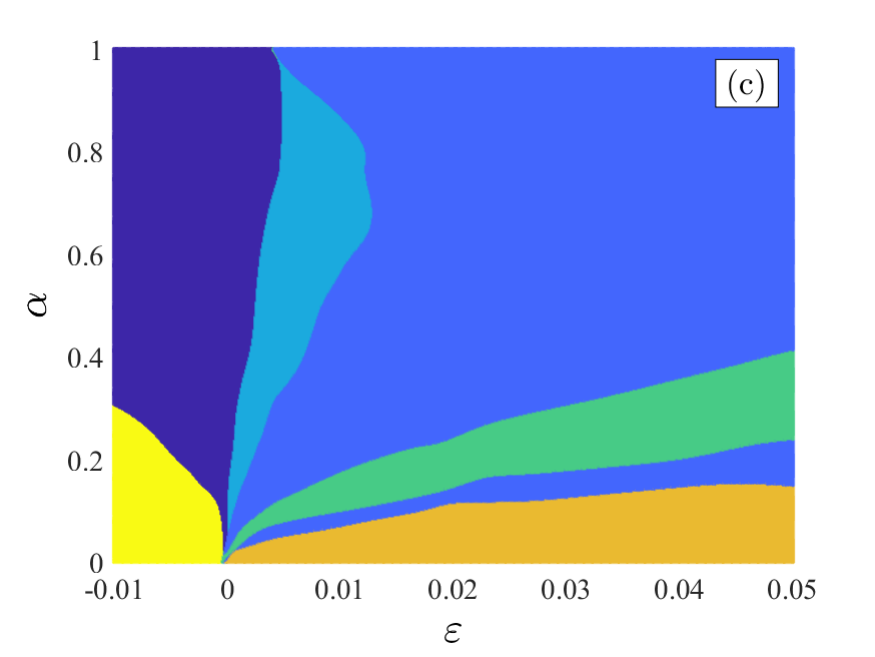}
    \caption{DiGCA phase classification accuracy with noise levels of 1\%(a), 5\%(b), and 10\%(c).}
    \label{fig:robustness}
\end{figure}

\begin{figure}[htbp]
\centering
\includegraphics[width=0.32\linewidth]{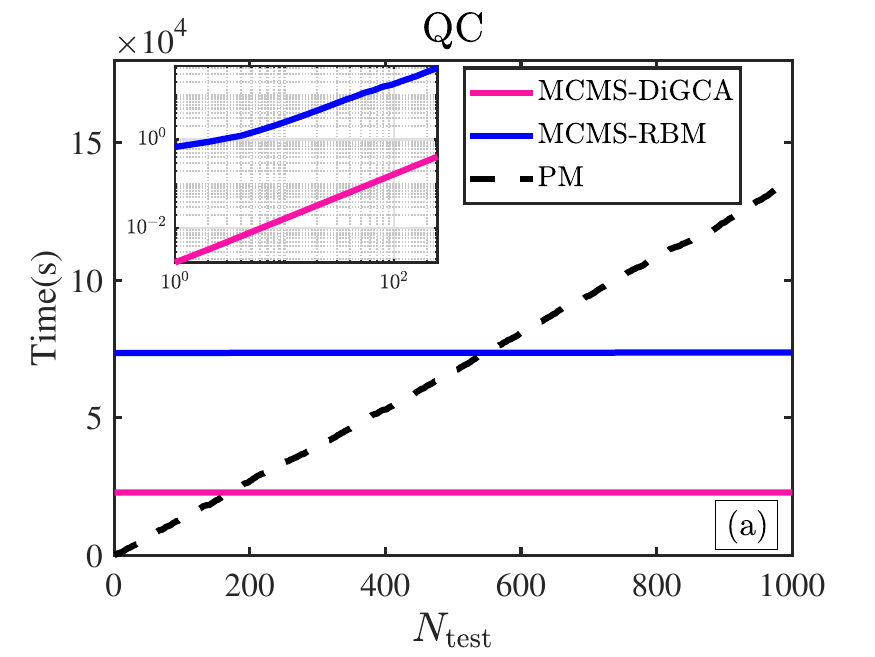}
\includegraphics[width=0.32\linewidth]{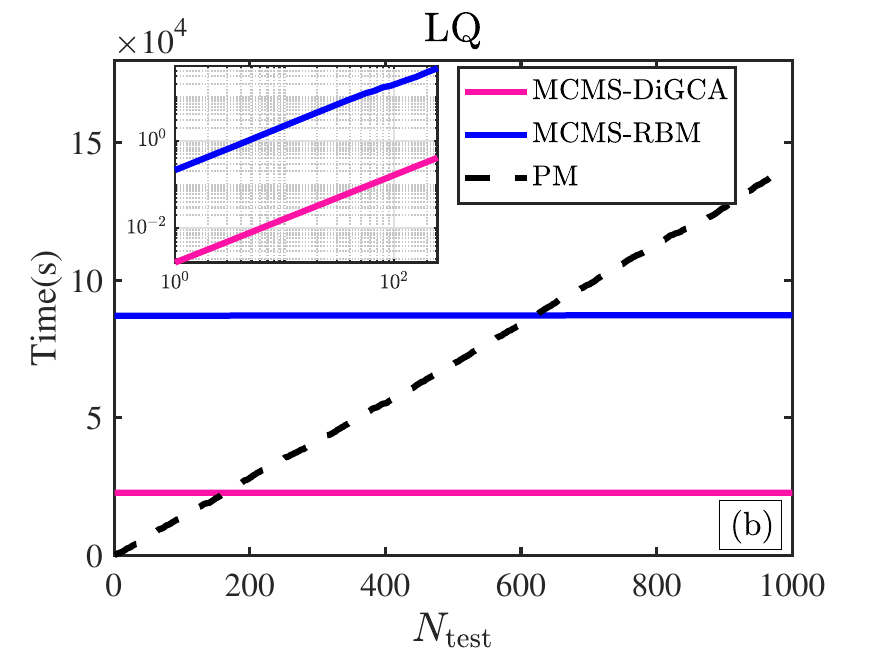}\\
\includegraphics[width=0.32\linewidth]{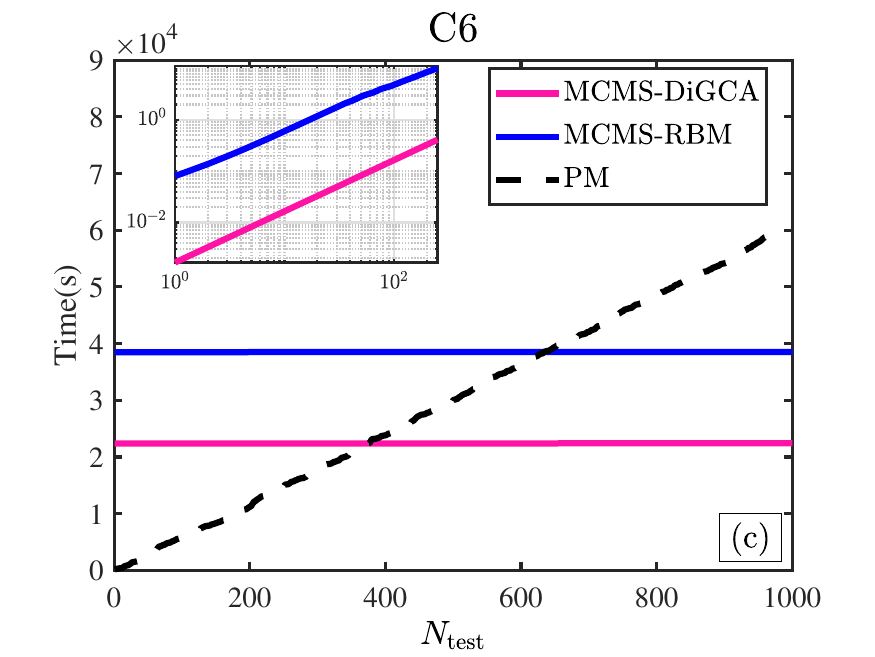}
\includegraphics[width=0.32\linewidth]{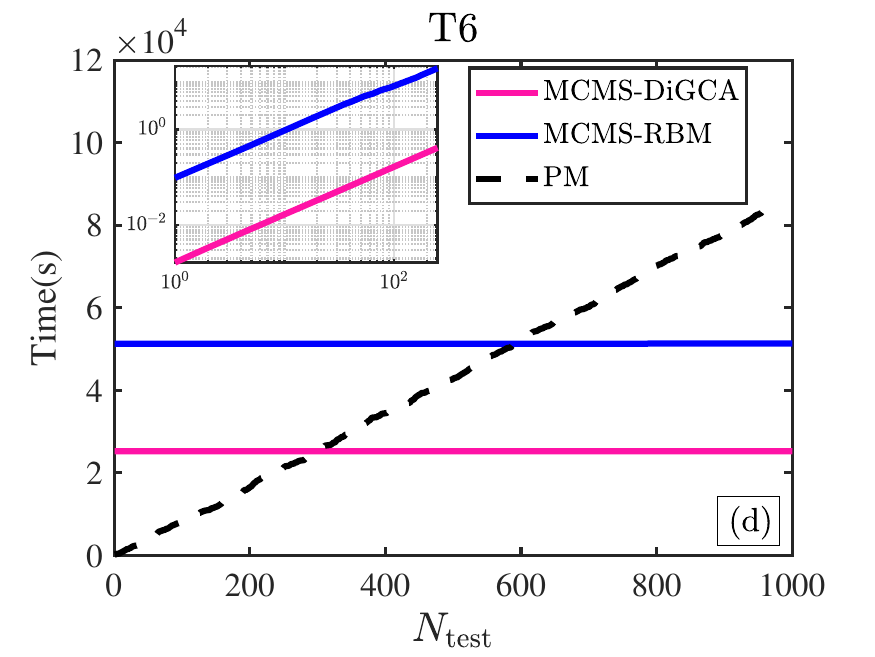}
\includegraphics[width=0.32\linewidth]{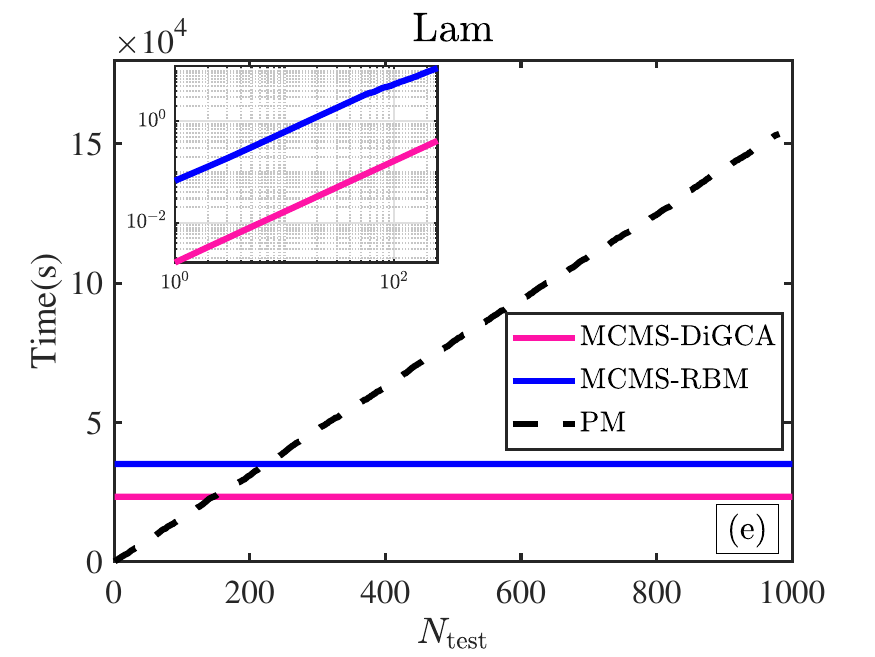}
\caption{Comparison of time performance between MCMS-RBM and MCMS-DiGCA.}
\label{fig:GAP_time}
\end{figure}

Our DiGCA-based classifier demonstrates superior robustness and computational efficiency compared with traditional methods while matching their accuracy. In Figure \ref{fig:robustness}, we show the phase diagram with increasingly higher levels of white noise added to the input data of DiGCA. The visually minimal change in the resulting diagram attests to the robustness of our approach. More importantly, the proposed method offers substantial computational advantages over existing approaches. In fact, we show in Figure \ref{fig:GAP_time} the cumulative runtime of the methods. It is clear that, compared to the MCMS-RBM, our DiGCA not only takes significantly less effort to train but also achieves approximately 100 times faster marginal computation in practice (see zoom-ins of each subplot). 
The combination of rapid convergence and efficient computation makes our approach particularly suitable for high-throughput materials discovery and experimental data analysis applications.

\section{Conclusion}
\label{sec:conclusion}
This paper proposes a Derivative-informed Graph Convolutional Autoencoder (DiGCA) with phase classification for the parametrized quasiperiodic LP model with two length scales. Featuring multiple components with each providing a graph convolutional autoencoder for one branch of the problem induced by one part of the parameter domain, the DiGCA with phase classification method serves as a generic framework for reduced order modeling of parametric problems whose solution has multiple states across the parameter domain. Compared to conventional methods, DiGCA with phase classification method achieves approximately two orders of magnitude speedup, robustness, and comparable level of accuracy, demonstrating its attractive performance for practical tasks of rapid generation of phase diagram.

\end{document}